# An integral equation for the identification of causal effects in nonlinear models

Wing Hung Wong
Departments of Statistics and Biomedical Data Science
Stanford University

**Abstract**

When the causal relationship between $X$ and $Y$ is specified by a structural equation, the causal effect of $X$ on $Y$ is the expected rate of change of $Y$ with respect to changes in $X$, when all other variables are kept fixed. This causal effect is not identifiable from the distribution of $(X,Y)$. We give conditions under which this causal effect is identified as the solution of an integral equation based on the distributions of $(X,Z)$ and $(Y,Z)$, where $Z$ is an instrumental variable.

**Introduction**

Suppose the causal relation between two real-valued randomly variables $X$ and $Y$ is specified by the structural equation $Y = f(X,U)$, where $U$ represents all other variables that may also affect $Y$. We assume $f(s,U)$ is smooth in $x$, and write $Y(x) = f(x,U)$, $Y^{(i)}(x) = \frac{\partial^i}{\partial x^i} f(x,U)$, i=1,2. We regard $\theta(x) = E(Y^{(1)}(x))$ as the causal effect of $X$ on $Y$. Discussion of this model and its relation to the potential outcome framework for causal inference was given in Wong (2021). Since $Y(x)$ and $Y^{(i)}(x)$ are not directly obtainable from $X$ and $Y$, $\theta(x)$ is not identifiable from the distribution $(X,Y)$ alone. The method of instrumental variable attempts to identify $\theta(x)$ from the joint distribution of $(X,Y,Z)$ where the instrumental variable $Z$ can affect $X$ through another equation $X = g(Z,V)$. However, identifiability results using instrumental variables are only available under very strong restrictions $f$ and $g$. These results and related literature had been reviewed in Wong (2021) and will not be repeated.

We consider the following nonlinear, nonparametric causal model

- $Y = f(X,U), \quad Y \in R, X \in R, U \in R^p, \; f$ is bounded and smooth in $x$     (1)
- $X = g(Z,V), \quad Z \in R^q, V \in R^r$     (2)
- $\sup_{x,z} p_z(x) < \infty$ where $p_z(\;)$ denotes the density function of $X(z)$     (3)
- $Z$ is independent of $(U,V)$     (4)

In (1), the condition that $f$ is bounded and smooth in $x$ means that $\sup_u |f(0,u)| < \infty$ and $\sup_u |\frac{\partial^i}{\partial x^i} f(x,u)| < m(x)$ for $i$=1, 2, where $m(\;)$ is a bounded and integrable function. Then, when $x \to \infty$, we have $Y(\infty) = \lim Y(x)$ exists and $E(Y(x)) \to E(Y(\infty))$. Similarly for $Y(-\infty)$. Also, $\theta(x) = E(Y^{(1)}(x))$ is a differentiable function and $\lim \theta(x)$=0 as $x \to \pm\infty$.

For nonlinear $f$ and $g$, the independence condition (4) is not sufficient for the identification of $\theta(x)$ from the distribution of $(X,Y,Z)$. Under the condition that changes in $Y$ caused by varying $X$ is uncorrelated to changes in $X$ caused by varying Z, conditional on $Z = z$, Wong (2021) showed that the distributions $(X,Z)$ and $(Y,Z)$ identify a related function $\psi(z) =$

$E(Y^{(i1)}(X)|Z = z)$. That paper also demonstrated by examples that sometimes the function $\theta(x)$ can be recovered from the function $\psi(z)$, but did not provide results on the direct identification of $\theta(x)$. To fill this gap, in this paper we derive an integral equation that can be used to identify $\theta(x)$ from the distributions of $(X, Z)$ and $(Y, Z)$.

**Result**

To formulate our main result, consider the following conditions

- $I(X(z) \leq x)$ is uncorrelated with $Y^{(1)}(x)$, for all $x, z$ (5)
- The set of distributions of $X|Z = z$, induced by varying z, is a complete set (6)

<u>Theorem</u>: If (1)-(6) hold, then $\theta$ is identifiable via the integral equation

$$\int K(z, x)\theta(x)dx = \mu(z) - \mu(0) \tag{7}$$

$$\text{where } K(z, x) = P(X \leq x|Z = 0) - P(X \leq x|Z = z)$$

$$\mu(z) = E(Y|Z = z)$$

<u>Proof</u>:

$$\mu(z) = E(Y|Z = z) = E(f(X, U)|Z = z) = E(f(g(z, V), U|Z = z)$$

$$= E(Y(X(z)))$$

$$= E \int \delta(x - X(z))Y(x)dx \tag{8}$$

Replacing the delta function $\delta(\ )$ by the $N(0, \sigma^2)$ density $\phi_\sigma(\ )$, we define

$$\mu_\sigma(z) = E \int \phi_\sigma(x - X(z))Y(x)dx \tag{9}$$

Since $Y(x) = Y(X(z)) + Y^{(1)}(X(z))(x - X(z)) + \frac{1}{2}Y^{(2)}(X(W))(x - X(z))^2$, where $W$ is an intermediate variable lying between $x$ and $X(z)$, hence

$$\mu_\sigma(z) = EY(X(z) + E[\frac{1}{2}Y^{(2)}(X(W)) \int \phi_\sigma((x - X(z))(x - X(z))^2 dx$$

$$= \mu(z) + \frac{\sigma^2}{2} \sup_x m(x)$$

Thus, $|\mu_\sigma(z) - \mu(z)| \leq c\sigma^2$ for some constant $c$ (10)

Next, we claim that

$$|E\left(\Phi\left(\frac{x-X(z)}{\sigma}\right)Y^{(1)}(x)\right) - P(X(z) \leq x)\theta(x)| \leq cm(x)\sqrt{\sigma} \text{ for some constant } c \tag{11}$$

Assuming (11) is true, we now analyze the integral in (9). Using integration by part, we have

$$\mu_\sigma(z) = E[Y(\infty) - \int \left(\Phi\left(\frac{x-X(z)}{\sigma}\right)Y^{(1)}(x)\right)dx]$$

$$= E(Y(\infty)) - \int P(X(z) \leq x)\theta(x)dx + r(z, \sigma)$$

where for some constant $c$, $|r(z, \sigma)| \leq c\sqrt{\sigma}$ for all small $\sigma$.

Thus $|(\mu_\sigma(z) - \mu_\sigma(0)) - \int [P(X(0) \leq x) - P(X(z) \leq x)]\theta(x)dx| \leq 2c\sqrt{\sigma}$ (12)

Taking the limit of (10) and (12) as σ→ 0, we have

$$\mu(z) - \mu(0) = \lim_{\sigma \to 0}(\mu_\sigma(z) - \mu_\sigma(0)) = \int [P(X(0) \le x) - P(X(z) \le x)]\theta(x)dx.$$

The desired equation (7) follows because $P(X(z) \le x) = P(g(z,V) \le x) = P(g(z,V) \le x|Z=z) = P(g(Z,V) \le x|Z=z) = P(X \le x|Z=z)$.

To prove the claim (11),

$$\left|E\left(\Phi\left(\tfrac{x-X(z)}{\sigma}\right)Y^{(1)}(x)\right) - P(X(z) \le x)\theta(x)\right|$$

$$= \left|E\left(\Phi\left(\tfrac{x-X(z)}{\sigma}\right)Y^{(1)}(x)\right) - E(I(X(z) \le x))E(Y^{(1)}(x))\right|$$

$$= \left|E\left(\Phi\left(\tfrac{x-X(z)}{\sigma}\right)Y^{(1)}(x)\right) - E(I(X(z) \le x)Y^{(1)}(x))\right| \quad \text{(by condition (5))}$$

$$\le m(x)\, E\left|\Phi\left(\tfrac{x-X(z)}{\sigma}\right) - I(X(z) \le x)\right|$$

$$\le m(x)\left[\Phi\left(-\tfrac{1}{\sqrt{\sigma}}\right) + 4(\sup_{x,z} p_z(x))\sqrt{\sigma}\,\right] \tag{13}$$

The last inequality (13) holds because $\left|\Phi\left(\tfrac{x-X(z)}{\sigma}\right) - I(X(z) \le x)\right|$ is bounded by 2 on $A(\sigma)$ and by $\Phi(-\tfrac{1}{\sqrt{\sigma}})$ on $A(\sigma)^C$, where $A(\sigma)$ is the event $\{|X(z) - x| \le \sqrt{\sigma}\}$.

Since both $K(z,x)$ and $\mu(z)$ in the integral equation (7) are determined by the distributions of $(X,Z)$ and $(Y,Z)$, it follows that $\theta$ is also determined if the solution to (7) is unique.

To establish uniqueness, let $a$ be a fixed constant, and define for any $\theta(\ )$, its anti-derivative $\lambda(x) = a - \int_x^\infty \theta(t)dt$. Suppose $\theta_1$ and $\theta_2$ are two solutions to (7) and $\lambda_1$ and $\lambda_2$ are the corresponding anti-derivatives, then

$$E(\lambda_1(X) - \lambda_2(X)|Z=z) = \int p_{X|Z}(x|z)\,(\lambda_1 - \lambda_2)(x)dx$$

$$= -\int P(X \le x|Z=z)(\theta_1 - \theta_2)(x)dx = -\int P(X \le x|Z=0)(\theta_1 - \theta_2)(x)dx.$$

Since the last expression does not depend on $z$, condition (6) implies $\lambda_1 = \lambda_2$, and therefore $\theta_1 = \theta_2$.

**Discussion**

Of the 6 conditions in the theorem, the first 3 are needed just set up the model and are not restrictive. On the other hand, conditions (4), (5), (6) each represents a significant constraint on the model. Condition (4) says that $Z$ is independent of all other causal variables that affect $X$ and $Y$. Together with (1) and (2), this means that the only way $Z$ can affect $Y$ causally is indirectly through its effect on $X$. This seems to be a natural condition on an instrumental variable. Condition (6) implies that the family of conditional distributions $P(X|Z=z)$ as $z$ varies, is a large family. This means that $Z$ has non-trivial relationship with $X$ in the sense that varying the value of $z$ leads to rich changes in the distribution of $X$. This is also a reasonable condition on an instrumental variable. This type of completeness condition was first introduced into causal inference by Imbens and Newey (2003). Finally, condition (5) requires $Y^{(1)}(x) = \frac{\partial f}{\partial x}(x, U)$ to be uncorrelated to $I(X(z) \leq x) = I(g(z, V) \leq x)$, which is a non-trivial condition not easy to interpret, but is needed to establish the relationship (7) between $\mu(z)$ and $\theta(x)$. Wong (2021) introduced a similar condition that requires $\frac{\partial f}{\partial x}(X, U)$ to be conditionally uncorrelated to $\frac{\partial g}{\partial z}(Z, V)$ given $Z = z$. However, under that condition one can only relate $\mu(z)$ to $\psi(z) = E(\frac{\partial f}{\partial x}(X(z), U))$ but not to $\theta(x) = E(\frac{\partial f}{\partial x}(x, U))$. In the general context of (1)-(4), we are not aware of alternative conditions that be used to relate $\mu(z)$ to $\theta(x)$.

Example: Suppose $Y = U_1 h(X) + U_2$, $X = g(Z, V)$, where $h( )$ is a smooth and bounded function in $x$. If $U_1$ is independent of $V$, then condition (5) is satisfied. In this example, the "subject-level" causal effect $Y^{(1)}(x)$ is assumed to be proportional to an nonlinear function $h(x)$, but heterogeneity is allowed by letting the proportionality constant depend on the subject. On the other hand, no restriction is imposed on the relation between $Z$ and $X$ beyond the completeness condition (6), which is not too restrictive. For example, (6) holds in the following cases (a) $g(z, v) = s(z + v)$ where $s( )$ is an invertible function and $V$ is a continuous random variable, (b) $g(z, v) = 1 + v_1 z + v_2 z^2$, $V_1$ and $V_2$ are independent random variables. This example demonstrated the usefulness of our result in a nonlinear, nonparametric model that allows heterogeneity in the causal effect of $X$ on $Y$ in different subjects.

The above proof of the theorem follows the way we discovered the integral equation originally, namely, start with the expression for $E(Y|Z = z)$, replace the delta function in the expression by the normal kernel and then integrate by part to obtain an expression involving $\theta( )$. Weijie Su (personal communication) suggests a second proof, which starts from the given $K(z, x)$ and then shows that the integral in (7) gives rise to $\mu(z) - \mu(0)$. His proof has the advantage that it does not require the existence of bounded second derivatives. See Su (2021, arXiv).

**Acknowledgment:** The author thanks Peng Ding and Weijie Su for helpful comments. This work was supported by NSF grants DMS1811920 and DMS1952386